\documentclass{amsart}
\usepackage{amscd,amssymb}

\theoremstyle{plain}
\newtheorem{thm}{Theorem}[section]
\newtheorem{prp}[thm]{Proposition}
\newtheorem{lem}[thm]{Lemma}
\newtheorem{cor}[thm]{Corollary}
\theoremstyle{remark}
\newtheorem{rem}[thm]{Remark}
\newcommand{\m}{\phantom{-}}
\newcommand{\x}{\negmedspace}
\newcommand{\Z}{\mathbb{Z}}
\newcommand{\R}{\mathbb{R}}
\newcommand{\C}{\mathbb{C}\mkern1mu}

\newcommand{\CP}{\mathbb{C\mkern1mu P}}
\newcommand{\Sph}{\mathbb{S}}
\DeclareMathOperator{\Sym}{Sym}
\DeclareMathOperator{\diag}{diag}

\newcommand{\1}{{\mathchoice
{\mathrm 1\mskip-4.2mu\mathrm l}{\mathrm 1\mskip-4.2mu\mathrm l}
{\mathrm 1\mskip-3.9mu\mathrm l}{\mathrm 1\mskip-4.0mu\mathrm l}}}
\newcommand{\bmat}{\begin{pmatrix}}
\newcommand{\emat}{\end{pmatrix}}
\newcommand{\bsmat}{\bigl(\begin{smallmatrix}}
\newcommand{\esmat}{\end{smallmatrix}\bigr)}
\newcommand{\Bsmat}{\Bigl(\begin{smallmatrix}}
\newcommand{\Esmat}{\end{smallmatrix}\Bigr)}
\newcommand{\bbsmat}{\biggl(\begin{smallmatrix}}
\newcommand{\eesmat}{\end{smallmatrix}\biggr)}
\newcommand{\BBsmat}{\Biggl(\begin{smallmatrix}}
\newcommand{\EEsmat}{\end{smallmatrix}\Biggr)}
\newcommand{\abs}[1]{\vert #1\vert}
\newcommand{\into}{\hookrightarrow}
\newcommand{\SO}{\mathrm{SO}}
\newcommand{\SU}{\mathrm{SU}}
\newcommand{\U}{\mathrm{U}}
\newcommand{\BU}{\mathrm{BU}}
\newcommand{\Syp}{\mathrm{Sp}}

\begin{document}
\title[Presentations of homotopy groups]{Presentations of the first
 homotopy groups\\ of the unitary groups}
\author{Thomas P\smash{\"u}ttmann}
\address{Ruhr-Universit\"at Bochum\\Fakult\"at f\"ur Mathematik\\
        D-44780 Bochum\\Germany}
\email{puttmann@math.ruhr-uni-bochum.de}
\author{A.~Rigas}
\address{IMECC\\UNICAMP, C.P.~6065\\ 13083-970 Campinas, SP\\Brazil}
\email{rigas@math.ime.unicamp.br}
\thanks{The joint work of the authors was supported by CNPq and
the International Bureau of the BMBF in the scope of the former
CNPq/GMD-agreement}


\begin{abstract}
We describe explicit presentations of all stable and the first nonstable
homotopy groups of the unitary groups. In particular, for each $n\ge 2$
we supply $n$ homotopic maps that each represent the $(n-1)!$-th
power of a suitable generator of $\pi_{2n}\SU(n)\approx\Z_{n!}$.
The product of these $n$ commuting maps is the constant map to
the identity matrix.
\end{abstract}

\maketitle

%
%

\section*{Introduction}
\label{s:intro}
The homotopy groups of compact Lie groups have been of continuous
interest since the discovery of homotopy groups at around 1935.
There is now a tremendous amount of computational tools
available and many groups have been determined.
On the other hand, intellectually and practically satisfying
presentations of these groups are only known in comparatively few cases.
Our goal in this paper is to describe such presentations
for the first homotopy groups of the unitary groups.

For the stable groups we mainly, but not entirely, review some known
results and procedures in an easily accessible and most explicit way.
Particular emphasis is given to the last stable groups $\pi_{2n-1}\U(n)$.
A highlight of this part is a strikingly simple formula for a minimal
embedding of $\Sph^5$ into $\SU(3)$ that represents a generator of
$\pi_5\SU(3)$ and has a natural interpretation
in terms of the complex cross product.

The main achievement of the paper concerns the first
nonstable homotopy groups $\pi_{2n}\SU(n)\approx \Z_{n!}$.
These groups played an important role in the first proofs
of the fact that the only parallelizable spheres are $\Sph^1$,
$\Sph^3$, and $\Sph^7$.
We figure in an elementary and explicit way how a suitable generator of
$\pi_{2n}\SU(n)$ becomes null-homotopic in the $n!$-th power.
Namely, we supply $n$ homotopic maps that each represent
the $(n-1)!$-th power of the generator. The product of these $n$
commuting maps is the constant map to the identity matrix.
The generators of $\pi_{2n}\SU(n)$ are then used to produce similar
presentations of the homotopy groups $\pi_{2n}\SU(n-1)$ with even $n$.

Finally, we answer the simple question whether a map $\Sph^r\to \U(n)$
is homotopic to its transposed or its complex conjugate for $r \le 2n+2$
with the help of the explicit maps described before.
As applications we obtain presentations of certain stable homotopy groups
of the symplectic groups and a structure theorem for certain nonstable
homotopy groups of the symmetric spaces $\SU(n)/\SO(n)$.

\smallskip

Throughout this paper we use the well-known fact that the homotopy group
$\pi_k(G)$ of a compact connected Lie group $G$ is isomorphic to the
group of free homotopy classes of maps $\Sph^k\to G$. Here, the product
between two free homotopy classes is given by multiplying the representing
maps value by value with the product of~$G$. We also often use the
elementary fact that the inclusions $\SU(n)\to\U(n)$ induce isomorphisms
between $\pi_r\SU(n)$ and $\pi_r\U(n)$ for $r \ge 2$.

\medskip

\section{The stable homotopy groups of the unitary groups}
\label{s:stable}
\subsection{Bott periodicity}
It has been known since around 1940 that the inclusion of $\U(n)$
into $\U(n+1)$ induces an isomorphism between the homotopy groups
$\pi_r\U(n)$ and $\pi_r \U(n+1)$ if $r < 2n$.
The homotopy groups in this range are called stable.
Their simple structure became visible at the end of the 50's by
Bott's famous periodicity theorem \cite{bottmorse}:
The stable groups $\pi_r\U$ are trivial if $r$ is even and
isomorphic to $\Z$ if $r$ is odd.
In fact, Bott constructed isomorphisms
\begin{gather*}
  \pi_r\U(n)\to \pi_{r+2}\SU(2n)
\end{gather*}
for $r < 2n$ and thus all stable groups are determined by
$\pi_1\U(1)\approx \Z$ and the trivial group $\pi_2\U(2)$.
The periodicity isomorphisms can be given in the following explicit way:
One assigns to a map $\theta: \Sph^{r} \to \U(n)$ the map
$\mathcal{B(\theta)}: \Sph^{r+2} \to \SU(2n)$
defined on the unit sphere in $\C\times\R^{r+1}$~by
\begin{gather*}
\label{e:period}
  \mathcal{B(\theta)}
  \bsmat w\\ x \esmat =
  \Bsmat \1 & 0\\ 0 & \theta(\hat x)\Esmat
  \Bsmat w\1 & -\abs{x}\1\\ \abs{x}\1 & \m\bar w\1\Esmat
  \Bsmat
    \1 & 0\\ 0 & \overline{\theta(\hat x)}^t
  \Esmat
  = \Bsmat
    w\1 & -\abs{x}\overline{\theta(\hat x)}^t\\
    \abs{x}\theta(\hat x) & \m\bar w\1
  \Esmat.
\end{gather*}
Here, $\hat x$ stands as an abbreviation for the unit vector
$\tfrac{x}{\abs{x}} \in \Sph^r$ and $\1$ denotes the
$n\times n$ identity matrix.
This assignment $\mathcal{B}$ provides the periodicity isomorphism.
We refer to \cite{eckmann},\cite{eckensign} and \cite{fomenko} for
(essentially) this form of $\mathcal{B}$.
In \cite{eckmann},\cite{eckensign} it is deduced by the relation to
Hurwitz-Radon matrices with the help of $K$-theory.
In \cite{fomenko}, Bott's original arguments \cite{bottmorse}
are turned into an explicit formula.

\subsection{Totally geodesic presentations}
Iterating the periodicity isomorphism $\mathcal{B}$ above
starting with the parametrization
\begin{gather*}
  \zeta_1 : \Sph^1 \to \U(1), \quad
  z \mapsto z
\end{gather*}
provides embeddings $\zeta_k: \Sph^{2k-1} \to \U(2^{k-1})$
that represent generators of the groups $\pi_{2k-1}\U(2^{k-1})$
and take values in $\SU(2^{k-1})$ if $k \ge 2$. For example,
\begin{gather*}
  \zeta_2 : \Sph^3 \to \SU(2), \quad
  \bsmat w\\ z\esmat \mapsto
  \bsmat w & -\bar z\\ z & \m \bar w\esmat
\end{gather*}
is the standard parametrization of $\SU(2)$
and
\begin{gather*}
  \zeta_3 : \Sph^5 \to \SU(4), \quad
  \Bsmat z_1\\ z_2\\ z_3\Esmat \mapsto
  \bbsmat z_1 & \m 0 & -\bar z_2 & -\bar z_3\\
    0 & \m z_1 & \m z_3 & -z_2\\
    z_2 & -\bar z_3 & \m \bar z_1 & \m 0\\
    z_3 & \m \bar z_2 & \m 0 & \m \bar z_1\eesmat.
\end{gather*}
The embeddings $\zeta_k$ are totally geodesic and $\R$-linear
in the sense that they extend to $\R$-linear maps from $\R^{2k}$
to the space of complex $2^{k-1}\times 2^{k-1}$ matrices.
By placing several copies of $\zeta_k$ or its inverse along the
diagonal in a sufficiently large square matrix one can realize
all elements of the homotopy group $\pi_{2k-1}\U$ by totally
geodesic, $\R$-linear embeddings.
For all these and additional facts we refer to \cite{eckmann},
\cite{eckensign}, and \cite{fomenko}.

\smallskip

We have just seen that the homotopy groups $\pi_{2k-1}\U(n)$ admit
very simple presentations if $k$ is very small compared to $n$.
The question we are now going to answer is how one can obtain
presentations of the last stable groups $\pi_{2n-1}\U(n)$
in the sequence $\pi_r\U(n)$ with fixed $n$.

\subsection{A deformation}
\label{s:deform}
Consider the subset of $\SU(n+1)$ that consists of matrices whose
lower right entry vanishes.
There is the following map from this subset to the group $\SU(n)$:
\begin{gather*}
  \bsmat
     A & b\\ \bar c^t\x & 0
  \esmat
  \mapsto A - b \bar c^t.
\end{gather*}
Here $A$ is an $n\times n$-matrix and $b$, $c$ are unit vectors in $\C^n$.
The map above can be obtained by the following deformation in $\SU(n+1)$:
\begin{gather}
\label{e:deform}
  \Bsmat
     A - b \bar c^t \sin t & b \cos t\\ \bar c^t \cos t \x & \sin t
  \Esmat.
\end{gather}
For $t = 0$ we get the initial matrix above and for $t=\tfrac{\pi}{2}$
we obtain the target matrix embedded in the upper left $n\times n$-block
of $\SU(n+1)$. This deformation exists analogously in $\SO(n+1)$
and $\Syp(n+1)$.

\subsection{Factorization of the periodicity isomorphism}
Because of stability the original periodicity isomorphism $\mathcal{B}$
admits the following factorization:
\begin{gather*}
  \pi_{r}\U(n) \to  \pi_{r+2}\SU(n+1) \to \pi_{r+2}\SU(2n).
\end{gather*}
In order to obtain this factorization explicitly,
we review essentially an algorithm of Lundell~\cite{lundell}.
There are, however, modifications in the details and we substitute
some of his arguments by the simple explicit deformation above.
The algorithm itself is very short:
In a first step one deforms the map $\mathcal{B}(\theta)$
with values in $\SU(2n)$ by multiplying the matrix
\begin{gather}
\label{e:exchange}
  \Bsmat \1 & 0 & 0\\ 0 & \cos t & -\sin t\\ 0 & \sin t& \m \cos t\Esmat
\end{gather}
from the left. For $t=\tfrac{\pi}{2}$ the lower two rows in each value
of $\mathcal{B}(\theta)$ are exchanged (one changes the sign). Hence,
the lower right entry of the resulting matrix valued map vanishes.
In a second step one now applies the deformation to $\SU(2n-1)$
described above. It is not complicated to check that these two
deformation steps can be iterated until the map takes values in $\SU(n+1)$.

\subsection{The last stable groups $\pi_{2n-1}\U(n)$}
Iterating the factorized version of the periodicity isomorphism
starting with the map $\eta_1 = \zeta_1$ above we obtain maps
$\eta_n$ that represent generators of the last stable groups
$\pi_{2n-1}\U(n)$ and take values in $\SU(n)$ for $n \ge 2$.
Note that $\eta_2 = \zeta_2$ still is the standard
parametrization of $\SU(2)$. In the case $n=3$ one
applies the deformation above to the generator $\zeta_3$ of $\pi_5\SU(4)$.
This yields the map
\begin{gather*}
  \eta_3 : \Sph^5 \to \SU(3),\quad
  \Bsmat z_1 \\ z_2 \\ z_3\Esmat \mapsto
  \bbsmat z_1 + \bar z_3 z_2 & -\bar z_3^2 & -\bar z_2 + \bar z_3\bar z_1\\
    z_2^2 & z_1 - z_2\bar z_3 & z_3 + z_2\bar z_1\\
    -z_3 + \bar z_1 z_2 & -\bar z_2 - \bar z_1\bar z_3 & \bar z_1^2\eesmat.
\end{gather*}
A map of this form was obtained by Chaves and Rigas \cite{chaves}
with a related but slightly more complicated approach.
With a few transformations we simplify the formula in a way that a
striking relation to the complex cross product appears and
its equivariance properties are revealed.
In fact, after multiplying $\eta_3$ by
$\Bsmat 0 & 0 & 1\\ 0 & 1 & 0\\ 1 & 0 & 0\Esmat$
from the left and by
$\Bsmat 0 & \m 1 & \m 0\\ 0 & \m 0 & -1\\ 1 & \m 0 & \m 0\Esmat$
from the right and after passing from $z_1$ and $z_3$ to $\bar z_1$
and $\bar z_3$ we obtain the map $\eta$ of the next section.
There, we introduce the resulting map directly by the complex cross product.
For $n \ge 4$ we do not know whether any of the maps $\eta_n$ or any
map homotopic to some $\eta_n$ has any nice geometric properties
or can be found by a more geometric construction.

\smallskip

We now point out a property of the generators of $\pi_{2n-1}\U(n)$
that will be the key for understanding the first nonstable groups
$\pi_{2n}\U(n)\approx\Z_{n!}$ in Section\,\ref{s:firstunstable}.
Given a map $\theta: \Sph^{2n-1}\to\U(n)$ we obtain a map
$p_j\circ\theta: \Sph^{2n-1}\to\Sph^{2n-1}$ using the projection
$p_j: \U(n)\to \Sph^{2n-1}$ that maps a matrix to its $j$-th column.
\begin{lem}
\label{l:fundproperty}
The assignment
\begin{gather*}
  \theta \mapsto \tfrac{1}{(n-1)!} \deg (p_j\circ\theta)
\end{gather*}
yields an isomorphism $\pi_{2n-1}\U(n)\to\Z$.
This isomorphism is independent of $j$.
\end{lem}

In other words, a map $\theta: \Sph^{2n-1}\to\U(n)$ represents a
generator of $\pi_{2n-1}\U(n)$ if and only if the composition
with the projection to some (and hence any) of the columns has
degree $\pm (n-1)!$ where the sign is independent of the column.

\begin{proof}
The first part follows immediately from the exact homotopy sequence
of the bundle $\U(n-1)\to\U(n)\to\Sph^{2n-1}$ using the fact that
$\pi_{2n-2}\U(n-1)\approx \Z_{(n-1)!}$ and stable homotopy groups.
In order to see that, say, $p_{n-1}$ and $p_n$ yield the same
isomorphism, we multiply the values of $\theta$ from the right by the
matrix in (\ref{e:exchange}) with $t=\tfrac{\pi}{2}$.
The resulting map $\theta'$ is homotopic to $\theta$
and we get $p_{n-1}\circ\theta' = p_n\circ\theta$.
\end{proof}

\medskip

\section{A minimal generator of $\pi_5\SU(3)$}
\label{s:minimal}
Given two vectors $z,w\in\C^3$ their cross product is defined to be
\begin{gather*}
  z\times w = \Bsmat \bar z_2\bar w_3 - \bar z_3\bar w_2\\
    \bar z_3\bar w_1 - \bar z_1\bar w_3\\
    \bar z_1\bar w_2 - \bar z_2\bar w_1\Esmat.
\end{gather*}
If $z$ and $w$ are unit vectors that are perpendicular with respect
to the standard hermitian inner product on $\C^3$ then $z\times w$
is the unique vector such that the matrix whose columns are
$z$, $w$, and $z\times w$ is contained in $\SU(3)$. Hence,
\begin{gather}
\label{e:equiv}
  (A\cdot z)\times (A\cdot w) = A\cdot (z\times w)
\end{gather}
for all $A\in \SU(3)$ and $z,w\in \C^3$. 

We can now define an embedding $\eta:\Sph^5 \to \SU(3)$ by setting
$\eta(z)\cdot \bar z = z$ and $\eta(z)\cdot \bar w = z\times w$
if $w$ is perpendicular to $z$. This map is obviously not null-homotopic.
For if it were homotopic to the constant map from $\Sph^5$ to the identity
in $\SU(3)$ then the map from $\Sph^5$ to itself given by complex conjugation
would be homotopic to the identity of $\Sph^5$, which is not true.
An explicit formula for $\eta$ is given as follows:
\begin{gather*}
  \eta(z) = z z^t + \Bsmat \m 0 & -\bar z_3 & \m \bar z_2\\
    \m \bar z_3 & \m 0 & -\bar z_1\\ -\bar z_2 & \m \bar z_1 &  \m 0\Esmat.
\end{gather*}

\begin{thm}
The embedding $\eta:\Sph^5 \to \SU(3)$ generates $\pi_5\SU(3)\approx\Z$.
\end{thm}

\begin{proof}
If we compose $\eta$ with the projection to any of the columns of $\SU(3)$
we obtain a map from $\Sph^5$ to itself with degree $2$. It follows from
Lemma\,\ref{l:fundproperty} that $\eta$ represents a generator of
$\pi_5\SU(3)$.
\end{proof}

It follows from property (\ref{e:equiv}) that $\eta$ is equivariant
with respect to the standard action of $\SU(3)$ on $\Sph^5\subset\C^3$
and the action of $\SU(3)$ on itself given as follows:
\begin{gather*}
  \SU(3) \times \SU(3) \to \SU(3), \quad
  (B,A) \mapsto BAB^t.
\end{gather*}
It is known that the orbit space of the latter action is a closed interval.
In order to give a more detailed description of the orbit structure
we use the geodesic
\begin{gather*}
  c(t) = \Bsmat 1 & 0 & 0\\ 0 & \cos t & -\sin t\\
    0 & \sin t& \m \cos t\Esmat.
\end{gather*}
The orbit through $c(0) = \1$ is diffeomorphic to the symmetric space $\SU(3)/\SO(3)$ and consists precisely of the symmetric matrices in $\SU(3)$.
It is easy to see that $c$ intersects this orbit perpendicularly
(and hence all orbits by Clairault's theorem that the velocity vectors
of a geodesic have a constant angle with a Killing field).
The orbits through $c(t)$ for $t \in \,]0,\tfrac{\pi}{2}[$ are
diffeomorphic to the seven-dimensional space $\SU(3)/\SO(2)$.
Finally, we have
\begin{gather*}
  c(\tfrac{\pi}{2})
  = \Bsmat 1 & \m 0 & \m 0\\ 0 & \m 0 & -1\\
    0 & \m 1 & \m 0\Esmat
  = \eta\Bsmat 1\\0\\0\Esmat.
\end{gather*}
Hence, $\eta$ parametrizes the isolated singular orbit through
$c(\tfrac{\pi}{2})$. Since isolated orbits are minimal submanifolds
(in the usual sense that they are critical points for the volume
functional, i.e., their mean curvature vanishes) \cite{lawson}, we get:

\begin{prp}
The embedding $\eta$ parametrizes a minimal submanifold of $\SU(3)$.
\end{prp}

We finally mention the following curiosity: Using the embedding $\eta$,
the Hopf fibration $\Sph^5 \to \CP^2$ can be extended to a simple
self-map of $\SU(3)$, namely, to the map $A\mapsto A\cdot \bar A$.
Indeed, if we multiply $\eta$ and $\bar\eta$ value by value
we obtain the map
\begin{gather*}
  \Sph^5 \to \SU(3),\quad
  z \mapsto 2 z\bar z^t - \1.
\end{gather*}
This map is the standard totally geodesic Cartan embedding of
$\CP^2$ into $\SU(3)$.
It follows from Theorem\,\ref{t:transposed} or by inspecting the orbit space
of the adjoint action of $\SU(3)$ that $\eta\cdot\bar\eta$ is null-homotopic.

\medskip

\section{The first nonstable homotopy groups $\pi_{2n}\SU(n)\approx\Z_{n!}$}
\label{s:firstunstable}
Bott \cite{bottloop} showed in 1958 that the image of $\pi_{2n}\BU$ in
$H_{2n}(\BU)$ is divisible by precisely $(n-1)!$. This refined the
previous result of Borel and Hirzebruch \cite{borel} that these classes
are divisible by $(n-1)!$ except for the prime $2$. As a consequence
of the refined version, the first nonstable homotopy groups $\pi_{2n}\U(n)$
of the unitary groups are isomorphic to the cyclic groups of order $n!$.
This result was used almost immediately by Kervaire \cite{kervparallel} and
Milnor \cite{bottmilnor} who independently gave the first proofs of the fact
that the only parallelizable spheres are $\Sph^1$, $\Sph^3$, and $\Sph^7$.
Generators of the groups $\pi_{2n}\U(n)$ are represented by the
characteristic maps of the bundles $\U(n+1) \to \Sph^{2n+1}$.
These maps were known explicitly several years before
Bott's result \cite{steenrod}.
We will deform them in a way that allows us to see
how they become null-homotopic in the $n!$-th power.

\smallskip

The group $\SU(n)$ acts transitively on the unit sphere $\Sph^{2n-1}$
in $\C^n$. The isotropy group of the $j$-th canonical basis vector in $\C^n$
is denoted by $\SU(n-1)_j$. It is the subgroup of $\SU(n)$ whose
$j$-th diagonal entry is $1$.
Natural diffeomorphisms between $\SU(n)/\SU(n-1)_j$ and $\Sph^{2n-1}$
are given by the projections $p_j: \SU(n)\to \Sph^{2n-1}$ that map
matrices to their $j$-th columns.
Now consider the maps
\begin{gather*}
  \phi_j : [0,\tfrac{2\pi}{n}] \times \SU(n)/\SU(n-1)_j \longrightarrow \SU(n)
\end{gather*}
given by
\begin{gather*}
  \phi_1(t,A) = A\cdot\diag(e^{i(n-1)t},e^{-it},\ldots,e^{-it})\cdot A^{-1}\\
  \vdots\\
  \phi_n(t,A) = A\cdot\diag(e^{-it},\ldots,e^{-it},e^{i(n-1)t})\cdot A^{-1}.
\end{gather*}
For $t=0$ and $t=\tfrac{2\pi}{n}$ the values of all $\phi_j$
are independent of $A\in\SU(n)$.
Hence the $\phi_j$ induce maps $\Sph^{2n}\to\SU(n)$.

\begin{lem}
All the maps $\phi_j$ above induce the same map $\phi: \Sph^{2n}\to\SU(n)$.
This map represents a generator of $\pi_{2n}\SU(n)$.
\end{lem}

\begin{proof}
Consider a matrix $A$ whose first column is given by $z\in\Sph^{2n-1}$.
Then
\begin{multline*}
  A\cdot\diag(e^{i(n-1)t},e^{-it},\ldots,e^{-it})\cdot A^{-1}\\
  = A \cdot e^{-it} \cdot
    \bigl(\1 + \diag(e^{int}-1,0,\ldots,0)\bigr)\cdot A^{-1}\\
  =  e^{-it}\bigl(\1 + z(e^{int}-1)\bar z^t\bigr).
\end{multline*}
For the other columns the computation is analogous and yields the
same result. We compose the map
\begin{gather*}
  \hat\phi: [0,\tfrac{2\pi}{n}] \times \Sph^{2n-1} \to \SU(n),\quad
  (t,z) \mapsto e^{-it}\bigl(\1 + z(e^{int}-1)\bar z^t\bigr)
\end{gather*}
with the inverse of the suspension
\begin{gather*}
  [0,\tfrac{2\pi}{n}] \times \Sph^{2n-1}
    \to \Sph^{2n} \subset \R\times\C^n, \quad
  (t,z) \mapsto
  \bigl( \tfrac{t n}{\pi} -1, 
    z \sqrt{\smash[b]{1-(\tfrac{t n}{\pi}-1)^2}}\bigr).
\end{gather*}
This yields the map
\begin{gather*}
  \phi: \Sph^{2n} \to \SU(n),\quad
  (y,z) \mapsto
    e^{-i\pi(y+1)/n}\cdot
    \bigl(\1 - \tfrac{z}{\abs{z}}(1+e^{i\pi y})\tfrac{\bar z^t}{\abs{z}}\bigr).
\end{gather*}
We can remove the factor in front of the paranthesis and obtain
a homotopic map with values in $\U(n)$.
Moreover, we can substitute the rational parametrization
$(\tfrac{1+iy}{1-iy})^2$ of the unit circle in $\C$
for the exponential parametrization $e^{i\pi y}$
without changing the homotopy class of $\phi$.
This leads to the map
\begin{gather*}
  \Sph^{2n} \to \U(n),\quad
  (y,z) \mapsto \1 - 2 z\tfrac{1}{(1-iy)^2}\bar z^t.
\end{gather*}
In Steenrod's book \cite{steenrod} it is proved that this map
represents the characteristic map of the bundle $\U(n+1)\to\Sph^{2n+1}$
and hence a generator of $\pi_{2n}\U(n)$.
\end{proof}

At first glance it might seem like one could multiply the $n$ maps
$\phi_1,\ldots,\phi_n$ value by value and the result is the constant map
to the identity.  This would imply that $\pi_{2n}\SU(n)$ is of order
at most $n$ contradicting $\pi_{2n}\SU(n)\approx\Z_{n!}$.
The reason why this does not work is that we are
not multiplying maps that have the same domain of definition,
since the isotropy groups $\SU(n-1)_j$ are different.
In order to get maps from the same domain of definition
$[0,\tfrac{2\pi}{n}]\times\Sph^{2n-1}$
one has to use the identifications between $\SU(n)/\SU(n-1)_j$
and $\Sph^{2n-1}$. But, as we saw, this always yields the same
map $\phi$ above and $\phi^n$ is evidently not the constant map
to the identity.

There is, however, a way to make the previous idea work.
The clue is to use any map $\eta: \Sph^{2n-1}\to \SU(n)$
that represents a generator of the stable group $\pi_{2n-1}\SU(n)$.
Such a map has the fundamental property that the composition
$p_j\circ\eta$ with the projection $p_j$ to the $j$-th matrix column
has degree $\pm (n-1)!$ where the sign is independent of the column
(see Lemma\,\ref{l:fundproperty}). We now obtain maps
\begin{gather*}
  \psi_j : [0,\tfrac{2\pi}{n}] \times \Sph^{2n-1} \to \SU(n)
\end{gather*}
by plugging $p_j\circ\eta$ into the second argument of $\hat\phi$,
i.e., by
\begin{gather*}
  \psi_1(t,z) = \phi_1(t,\eta(z))
  = \eta(z) \cdot \diag(e^{i(n-1)t},e^{-it},\ldots,e^{-it})
    \cdot \eta(z)^{-1}\\
  \vdots\\
  \psi_n(t,z) = \phi_n(t,\eta(z))
  = \eta(z) \cdot \diag(e^{-it},\ldots,e^{-it},e^{i(n-1)t})
   \cdot \eta(z)^{-1}.
\end{gather*}
The following is now evident and shows us explicitely how the
$n!$-th power of a generator of $\pi_{2n}\SU(n)$ is null-homotopic.

\begin{thm}
\label{t:main}
The maps $\psi_j$ induce maps $\Sph^{2n}\to\SU(n)$ that
represent $(n-1)!$ times the same generator of $\pi_{2n}\SU(n)$.
The maps $\psi_j$ commute mutually and their product
$\psi_1\cdot\ldots\cdot\psi_n$ is the constant map to the identity.
\end{thm}

Explicit homotopies between the maps $\psi_j$ are easily given.
The formula
\begin{gather*}
  \eta(z) \cdot
  \Bsmat \cos s & -\sin s & \m 0\\
    \sin s & \m \cos s & \m 0\\ 0 & \m 0 & \m \1 \Esmat
  \cdot \diag(e^{i(n-1)t},e^{-it},\ldots,e^{-it}) \cdot
  \Bsmat \m \cos s & \m \sin s & \m 0\\
    -\sin s & \m \cos s & \m 0\\ 0 & \m 0 & \m \1 \Esmat
  \cdot \eta(z)^{-1},
\end{gather*}
for example, yields the map $\psi_1$ for $s=0$ and the map $\psi_2$
for $s=\tfrac{\pi}{2}$.

\begin{rem}
Theorem\,\ref{t:main} and Lemma\,\ref{l:fundproperty} together provide
inductively an elementary proof for the fact that $\pi_{2n}\SU(n)$
is a cyclic group whose order devides $n!$. In order to show
that $n!$ devides the order of $\pi_{2n}\SU(n)$, however, cohomological
arguments like those in \cite{bottloop} seem to be inevitable.
\end{rem}

\begin{rem}
The map $\hat\phi$ factors through a map $\check\phi$ defined on
$[0,\tfrac{2\pi}{n}] \times \CP^{n-1}$. The $\CP^{n-1}$ can be
considered to represent the space of shortest curves from
the identity matrix $\1$ to the matrix $e^{-2\pi i/n}\cdot\1$ in the
center of $\SU(n)$. The map $\check\phi$ appears in Bott's papers \cite{bottloop},\cite{bottmorse} frequently, but not with the meaning
that it provides a generator of the group $\pi_{2n}\SU(n)$.
\end{rem}

\medskip

\section{The homotopy groups $\pi_{2n}\SU(n-1)$}
\label{s:secondunstable}
The homotopy groups $\pi_{2n}\SU(n-1)$ were first computed by Kervaire
\cite{kervaire}. The following fact is central for the computation:
Given a generator $\phi$ of $\pi_{2n}\SU(n)$ the composition
$p_j\circ\phi:\Sph^{2n}\to \Sph^{2n-1}$ with the projection to the
$j$-th column is null-homotopic if $n$ is odd and homotopic to the
$(2n-3)$-rd suspension of the Hopf fibration $\Sph^3\to\Sph^2$ if
$n$ is even (see \cite{kervaire}, \cite{steenrod}).
With this fact one easily deduces from the exact homotopy sequence
of the bundle $\SU(n)\to\Sph^{2n-1}$ that $\pi_{2n-1}\SU(n-1)$ is
trivial if $n$ is even and isomorphic to $\Z_2$ if $n$ is odd
and that
\begin{gather*}
  \pi_{2n}\SU(n-1) \approx
  \begin{cases}
    \Z_{n!/2} & \text{if $n$ is even,}\\
    \Z_2 \oplus \Z_{n!} & \text{if $n$ is odd.}
  \end{cases}
\end{gather*}
In other words, if $n$ is odd, $\phi$ is homotopic to a map with values in
$\SU(n-1)$ (it is, however, not very easy to write this homotopy down
explicitly). If $n=2m$ is even, $\phi$ cannot be deformed to a map with
values in $\SU(n-1)$, but $\phi^2$ can, and the resulting map represents
a generator of $\pi_{4m}\SU(2m-1)$.
We will now describe this deformation explicitly by reducing the
equivariance group of $\phi$ from $\SU(2m)$ to~$\Syp(m)$.

\smallskip

The symplectic group $\Syp(m)$ can be regarded as the subgroup of matrices
$A\in\SU(2m)$ with $A^tJA = J$. Here, $J\in\SU(2m)$ is the matrix
whose diagonal $2\times 2$-blocks are $\bsmat 0 & -1\\1 & \m 0\esmat$.
A matrix $A\in\SU(2m)$ with columns $v_1,\ldots,v_{2m}$
belongs to $\Syp(m)$ if and only if $v_{2k} = J\cdot \bar v_{2k-1}$
for all $k = 1,\ldots,m$.
The group $\Syp(m)$ acts transitively on the unit sphere $\Sph^{4m-1}$
in $\C^{2m}$. The isotropy groups of the first and the second canonical
basis vector in $\C^{2m}$ are the same, namely, the subgroup $\Syp(m-1)_1$
of matrices in $\Syp(m)$ whose first and second diagonal entry is $1$.
Now we consider the maps $\phi_1$ and $\phi_2$ of the previous section
and restrict the second argument of these maps to symplectic matrices.
This way we obtain maps
\begin{gather*}
   \phi_1',\phi_2' : [0,\tfrac{\pi}{m}] \times \Syp(m)/\Syp(m-1)_1 \to \SU(2m)
\end{gather*}
with the same domain of definition. Both these maps still induce the
generator $\phi : \Sph^{4m}\to\SU(2m)$ of $\pi_{4m}\SU(2m)$
given in the previous section.
Their product is the map
\begin{multline*}
   \phi_{12}' : [0,\tfrac{\pi}{m}] \times \Syp(m)/\Syp(m-1)_1 \to \SU(2m),\\
  (t,A) \mapsto
    A \cdot \diag(e^{i(2m-2)t},e^{i(2m-2)t},e^{-2it},\ldots,e^{-2it}) \cdot A^{-1}.
\end{multline*}

\begin{lem}
The map $\phi_{12}'$ and the analogously defined maps
$\phi_{34}',\ldots,\phi_{2m-1,2m}'$ all induce the same map
$\phi^{(2)} : \Sph^{4m} \to \SU(2m)$ which represents
twice a generator of $\pi_{4m}\SU(2m)$.
The deformation to $\SU(2m-1)$ of Section\,\ref{s:deform}
can be applied to $\phi^{(2)}$ and the deformed map represents
a generator of $\pi_{4m}\SU(2m-1)$.
\end{lem}

\begin{proof}
Since $\phi_{12}'$ is the product of $\phi_1'$ and $\phi_2'$
it is evident that $\phi^{(2)}$ represents twice a generator
of $\pi_{4m}\SU(2m)$.
Computations analogous to that of the previous section show the following:
If the first column of the matrix $A$ is the vector $z\in\C^{2m}$ then
\begin{gather*}
  \phi_{12}'(t,A) = e^{-2it} \bigl(\1 + (e^{2m i t}-1)
   (z\bar z^t - J\bar z z^t J)\bigr),
\end{gather*}
and the map $\phi_{12}'$ induces the map
\begin{multline*}
  \phi^{(2)}: \Sph^{4m} \to \SU(2m),\\
  (y,z) \mapsto
    e^{-i\pi(y+1)/m}\cdot
    \bigl(\1 - (1+e^{i\pi y})(\tfrac{z}{\abs{z}}\tfrac{\bar z^t}{\abs{z}}
     - J \tfrac{\bar z}{\abs{z}} \tfrac{z^t}{\abs{z}} J)\bigr).
\end{multline*}
The $(2,1)$-entry in the values of $\phi^{(2)}$ is always zero.
Hence, after multiplying $\phi^{(2)}$ from the left and the right with
suitable permutation matrices, the lower right entry vanishes and the
deformation of Section\,\ref{s:deform} can be applied.
\end{proof}

Analogously to the previous section we can plug a generator
of $\pi_{4m-1}\Syp(m)\approx\Z$ into the second argument of
the maps $\phi_{12}',\ldots,\phi_{2m-1,2m}'$. The resulting
maps will be denoted by $\psi_{12}',\ldots,\psi_{2m-1,2m}'$.
Like $\phi_{12}'$ they can all be deformed to $\SU(2m-1)$
with the explicit deformation of Section\,\ref{s:deform}.

\begin{prp}
\label{p:secondunstable}
The $m$ maps $\psi_{2k-1,2k}'$ induce maps
$\Sph^{4m}\to \SU(2m-1)$ that represent $(2m-1)!$ times a
generator of $\pi_{4m}\SU(2m-1)$ if $m$ is odd and
$2\cdot(2m-1)!$ times a generator if $m$ is even.
They commute and their product is the constant map
to the identity.
\end{prp}

\begin{proof}
It follows from the exact homotopy sequence of the bundle
$\Syp(m)\to\Sph^{4m-1}$ that the composition of a generator
of $\pi_{4m-1}\Syp(m-1)$ with the projection to any of the
columns of $\Syp(m-1)$ yields a self-map of $\Sph^{4m-1}$
whose degree is the order of the cyclic group $\pi_{4m-2}\Syp(m-1)$.
Kervaire \cite{kervpont} first showed that the order of this group is 
$(2m-1)!$ if $m$ is odd and $2\cdot(2m-1)!$ if $m$ is even. 
\end{proof}

\begin{rem}
The map $\phi^{(2)}$ is homotopic to the map
\begin{gather*}
  \Sph^{4m} \to \U(2m),\quad
  (y,z) \mapsto
  \1 - \tfrac{2}{(1-iy)^2} (z\bar z^t - J\bar z z^t J)
\end{gather*}
with values in $\U(2m)$.
\end{rem}

\begin{rem}
The map $\phi_{12}'$ factors through a map defined on
\begin{gather*}
  [0,\tfrac{\pi}{m}] \times \Syp(m)/(\Syp(m-1)\times\Syp(1)).
\end{gather*}
In the case $m=1$ the factor on the right is trivial and the
maps $\phi_{12}'$ and $\phi^{(2)}$ are the constant maps to the identity.
In the case $m=2$ the factor on the right is diffeomorphic to $\Sph^4$
and hence $\phi_{12}'$ induces a map $\Sph^5\to \SU(4)$.
It is not difficult to see that this is the map $\zeta_3$ from
Section\,\ref{s:stable}. This means that a generator of $\pi_8\SU(3)$
is given by composing the first suspension of the Hopf fibration
$\Sph^7\to\Sph^4$ with the generator of $\pi_5\SU(3)$ described
in Section\,\ref{s:minimal}.
\end{rem}

\medskip

\section{Symmetric maps into the unitary groups and homotopy groups\\
of $\Syp(n)$ and $\SU(n)/\SO(n)$}
The Cartan embedding of the symmetric space $\SU(n)/\SO(n)$
into the Lie group $\SU(n)$ is the map
\begin{gather*}
  \mathcal{C} : \SU(n)/\SO(n) \to \SU(n), \quad
  A\cdot \SO(n) \mapsto A\cdot A^t.
\end{gather*}
The image of this map is precisely the space of symmetric
matrices in $\SU(n)$. We combine this fact  with the explicit form
of the Bott periodicity isomorphism given in Section\,\ref{s:stable}.
This combination provides first maps that represent non-trivial
elements of certain homotopy groups of the symplectic groups and
second a structure theorem for certain nonstable homotopy groups
of $\SU(n)/\SO(n)$.

\smallskip

We first use the generators for the homotopy groups of the unitary groups
given in Section\,\ref{s:stable} and Section\,\ref{s:firstunstable}
to derive the following statement:
\begin{thm}
\label{t:transposed}
Any map $\Sph^{2k-1}\to\U(n)$ with $k\le n$ is
homotopic to its transposed if $k$ is odd and homotopic to its
complex conjugate if $k$ is even. Any map from $\Sph^{2n}$ to $\U(n)$
or $\U(n-1)$ is homotopic to its transposed if $n$ is even and homotopic
to its complex conjugate if $n$ is odd.
\end{thm}

\begin{proof}
Since the product in the homotopy groups $\pi_r\U(n)$ is induced
by the matrix product of $\U(n)$ it suffices to proof the properties
for maps that represent generators of $\pi_r(\U(n))$.
In the stable range $r=2k-1 \le n$ we consider the generators $\zeta_k$
of $\pi_{2k-1}\U(2^{k-1})$ that were given in Section\,\ref{s:stable}.
Inductively, we see that each map $\zeta_k$ satisfies
$\overline{\zeta_k(z)} = \zeta_k(\bar z)$ for all $z\in \Sph^{2k-1}$.
If $k$ is even the map $\Sph^{2k-1}\to \Sph^{2k-1}$, $z \mapsto \bar z$ is
an orthogonal transformation of the real vector space $\C^k\approx\R^{2k}$
with determinant $1$ and therefore homotopic to the identity.
Hence, $z \mapsto \zeta_k(z)$ is homotopic to
$z \mapsto \zeta_k(\bar z)= \bar\zeta_k(z)$.
It is now easy to see that $\zeta_{k+1}= \mathcal{B}(\zeta_k)$
is homotopic to its transposed. One has to use the fact that
$\zeta_k$ is homotopic to $-\zeta_k$.

Similarly, we see that the generator $\phi:\Sph^{2n}\to\U(n)$ of
$\pi_{2n}\U(n)$ given in Section\,\ref{s:firstunstable} is homotopic
to $\phi^t$ if $n$ is even and homotopic to $\bar \phi$ if $n$ is odd.
A generator of the homotopy group $\pi_{2n}\U(n-1)$
is homotopic to a generator of $\pi_{2n}\U(n)$ or to twice
a generator or of order $2$ (see Section\,\ref{s:secondunstable}).
\end{proof}

Note that in $\SU(2)$ complex conjugation is an inner automorphism
and thus any non-trivial map from any sphere $\Sph^m$ to $\SU(2)$
is homotopic to its conjugate and therefore only homotopic to its
transposed if it is of order $2$. From the higher homotopy groups
of $\SU(2)\approx\Sph^3$ it is now clear that there exist many maps
$\Sph^{2k-1}\to \SU(2)$ with odd $k$ that are not homotopic to their
transposed.

\smallskip

We now state an elementary property of the explicit form of the
Bott periodicity isomorphism given in Section\,\ref{s:stable}.
\begin{lem}
Let $\theta$ be a symmetric map $\Sph^{2k-1}\to\U(n)$, i.e.,
$\theta(z)^t = \theta(z)$ for all $z\in\Sph^{2k-1}$.
Then $\mathcal{B}(\theta)$ takes values in $\Syp(n)\subset\SU(2n)$.
\end{lem}
\begin{proof}
Since $\theta$ is symmetric, the image $\mathcal{B}(\theta)$
is a matrix of the form
$\bsmat A & -\bar B\\B & \m \bar A\esmat$.
Matrices of this form build a standard $\Syp(n)$ in $\SU(2n)$.
\end{proof}

\begin{cor}
Any symmetric map $\theta:\Sph^{2k-1}\to\U(n)$ with even $k\le n$
is null-homotopic.
\end{cor}

\begin{proof}
If $\theta$ would represent a non-trivial element in the stable group
$\pi_{2k-1}\U(n)$, then $\mathcal{B}(\theta)$ would represent a
non-trivial element in $\pi_{2k+1}\SU(2n)\approx\Z$.
But $\mathcal{B}(\theta)$ takes values in $\Syp(n)$ and the stable group
$\pi_{2k+1}\Syp(n)$ is trivial or isomorphic to $\Z_2$ if $k$ is even.
\end{proof}

\begin{prp}
If $\theta$ is a generator of $\pi_{2k-1}\U(n)$ with $k=2m+1 \le n$
then $\mathcal{B}(\theta\cdot\theta^t)$ represents a generator
of $\pi_{2k+1}\Syp(n)\approx\Z$ if $m$ is odd and
twice a generator if $m$ is even.
\end{prp}

\begin{proof}
Because of Theorem\,\ref{t:transposed}, $\theta\cdot\theta^t$ represents
in $\pi_{2k-1}\U(n)$ twice the generator given by $\theta$.
Correspondingly, $\mathcal{B}(\theta\cdot\theta^t)$ represents twice a
generator of $\pi_{2k+1}\SU(2n)$. Since $\theta\cdot\theta^t$
is symmetric, $\mathcal{B}(\theta\cdot\theta^t)$ falls into $\Syp(n)$.
We now inspect part of the exact homotopy sequence of the
bundle that belongs to the homogeneous space $\SU(2n)/\Syp(n)$:
\begin{gather*}
  \pi_{2k+1}\Syp(n) \to \pi_{2k+1}\SU(2n) \to
  \pi_{2k+1}\bigl(\SU(2n)/\Syp(n)\bigr) \to \pi_{2k}\Syp(n).
\end{gather*}
All the homotopy groups involved here are stable, the first two
isomorphic to $\Z$, the last one trivial, and
$\pi_{2k+1}\bigl(\SU(2n)/\Syp(n)\bigr)$ trivial if $m$ is even and
isomorphic to $\Z_2$ if $m$ is odd.
\end{proof}

\begin{cor}
A generator of the stable group $\pi_{2k-1}\U(n)$ with
$k = 4l+3 \le n$ cannot be represented by a symmetric map.
\end{cor}

We will now apply the statements above to determine the structure of
the semistable homotopy groups $\pi_{8l+5}\bigl(\SU(n)/\SO(n)\bigr)$
for $l \le \tfrac{n-3}{4}$.
Given a map $\theta:\Sph^{2k-1}\to\SU(n)$ with odd $k$,
the composition
\begin{gather}
\label{e:factor}
  \Sph^{2k-1} \stackrel{\theta}{\longrightarrow} \SU(n) \to \SU(n)/\SO(n)
  \stackrel{\mathcal{C}}{\longrightarrow} \SU(n)\cap\Sym(n,\C) \into \SU(n)
\end{gather}
yields the symmetric map $\theta\cdot\theta^t$ which represents twice
the element in $\pi_{2k-1}\SU(n)$ that is represented by $\theta$.
Hence, if $\theta$ represents a generator of $\pi_{2k-1}\SU(n)$ then $\theta\cdot\theta^t$ represents a generator or twice a generator of a
$\Z$-factor in
\begin{gather*}
  \pi_{2k-1}\bigl(\SU(n)\cap\Sym(n,\C)\bigr)
    \approx \pi_{2k-1}\bigl(\SU(n)/\SO(n)\bigr).
\end{gather*}
This ambiguity remains in the case $k=4l+1\le n$ as we shall see below.
However, if $k=4l+3\le n$ then $\theta\cdot\theta^t$ cannot represent
twice a generator because of the previous corollary.

\begin{thm}
If $3\le k = 4l+3\le n$ then a generator of $\pi_{2k-1}\SU(n)$ projects
to a generator of a $\Z$-factor in $\pi_{2k-1}\bigl(\SU(n)/\SO(n)\bigr)$.
Consequently we have
\begin{align*}
  \pi_{8l+5}\bigl(\SU(n)/\SO(n)\bigr) &\approx \Z \oplus \pi_{8l+4}\SO(n)
  \text{ and}\\
  \pi_{8l+6}\bigl(\SU(n)/\SO(n)\bigr) &\approx \pi_{8l+5}\SO(n).
\end{align*}
\end{thm}

\begin{proof}
The first part follows from the factorization (\ref{e:factor}) of the map
$\theta\cdot\theta^t$, the second part from the first and the relevant
segment of the exact homotopy sequence of the bundle
$\SO(n)\to\SU(n)\to\SU(n)/\SO(n)$.
\end{proof}

This statement was obtained before by Kachi
(see \cite{kachi}, Proposition\,3.5) for $7 \le k = 4l+3 \le n-1$.
Kachi's proof is based on computations of Kervaire \cite{kervaire}.
These, in turn, involve certain homotopy groups of the Stiefel manifolds
that were determined by Paechter \cite{paechter}.
Our proof, on the other hand, requires just the knowledge
of stable homotopy groups.

The simplest example where our statement provides information
is the homotopy group $\pi_5\bigl(\SU(3)/\SO(3)\bigr)$.
The exact homotopy sequence leaves the two choices $\Z$ and $\Z\oplus\Z_2$.
Our argument above shows that $\pi_5\bigl(\SU(3)/\SO(3)\bigr)$
is isomorphic to the second group.

\begin{prp}
If $5\le k = 4l+1\le n$ then a generator of $\pi_{2k-1}\SU(n)$
can only project to a generator or to twice a generator of a
$\Z$-factor in $\pi_{2k-1}\bigl(\SU(n)/\SO(n)\bigr)$.
In the first case we have
\begin{gather*}
  \pi_{8l+1}\bigl(\SU(n)/\SO(n)\bigr) \approx \Z \oplus \pi_{8l}\SO(n).
\end{gather*}
In the second case, $\pi_{8l}\SO(n)$ is isomorphic to a direct sum
$G \oplus \Z_2$ such that
\begin{gather*}
  \pi_{8l+1}\bigl(\SU(n)/\SO(n)\bigr) \approx \Z \oplus G.
\end{gather*}
In any of the two cases we have
\begin{gather*}
  \pi_{8l+2}\bigl(\SU(n)/\SO(n)\bigr) \approx \pi_{8l+1}\SO(n).
\end{gather*}
\end{prp}

In the stable range obviously the second alternative holds.
On the other hand, the group $\pi_8\SO(5)$ (which occurs in
the first case covered by the proposition) is trivial.
For information on $\pi_{n+r}\bigl(\SU(n)/\SO(n)\bigr)$
for $n\ge 8$ and $r\le 5$ we refer to~\cite{kachi}.

\bigskip

\appendix

\section*{Appendix. The first homotopy groups of the unitary groups}
For the convenience of the reader we provide a table of the very first
homotopy groups of the unitary groups. Larger tables can be found in
\cite{lundtable}. The black line in the table indicates the border
between the stable and the non stable groups. A $t$ or $c$ below a group
$\pi_{r}\U(n)$ indicates that any map $\Sph^r\to\U(n)$ is homotopic
to its transposed or to its complex conjugate, respectively.

\begin{table}
\begin{center}
\begin{tabular}{|c|cccccc|}
\hline
$r\;\backslash\; n$ & $\quad 1\quad$ & $\quad 2\quad$ & $\quad 3\quad$ & $\quad 4\quad$ & $\quad 5\quad$ & $\quad 6\quad$\\
\hline
$1$ & \multicolumn{1}{|c}{$\underset{t}{\Z}$} & $\Z$ & $\Z$ & $\Z$ & $\Z$ & $\Z$\\
\cline{2-2}
$2$ & $0$ & \multicolumn{1}{|c}{$0$} & $0$ & $0$ & $0$ & $0$\\
$3$ & $0$ & \multicolumn{1}{|c}{$\underset{c}{\Z}$} & $\Z$ & $\Z$ & $\Z$ & $\Z$\\
\cline{3-3}
$4$ & $0$ & $\underset{t\, c}{\smash[b]{\Z_2}}$ & \multicolumn{1}{|c}{$0$} & $0$ & $0$ & $0$\\
$5$ & $0$ & $\underset{t\, c}{\smash[b]{\Z_2}}$ & \multicolumn{1}{|c}{$\underset{t}{\Z}$} & $\Z$ & $\Z$ & $\Z$\\
\cline{4-4}
$6$ & $0$ & $\underset{\x\x\x c}{\smash[b]{\Z_{12}}}$ & $\underset{c}{\smash[b]{\Z_{6}}}$ & \multicolumn{1}{|c}{$0$} & $0$ & $0$\\
$7$ & $0$ & $\underset{t\, c}{\smash[b]{\Z_2}}$ & $0$ & \multicolumn{1}{|c}{$\underset{c}{\Z}$} & $\Z$ & $\Z$\\
\cline{5-5}
$8$ & $0$ & $\underset{t\, c}{\smash[b]{\Z_2}}$ & $\underset{\x\x\x t}{\smash[b]{\Z_{12}}}$ & $\underset{\x\x\x t}{\smash[b]{\Z_{24}}}$ & \multicolumn{1}{|c}{$0$} & $0$\\
$9$ & $0$ & $\underset{c}{\smash[b]{\Z_3}}$ & $\underset{c}{\smash[b]{\Z_3}}$ & $\underset{t\, c}{\smash[b]{\Z_2}}$ & \multicolumn{1}{|c}{$\underset{t}{\Z}$} & $\Z$\\
\cline{6-6}
$10$ & $0$ & $\underset{\x\x\x c}{\smash[b]{\Z_{15}}}$ & $\underset{\x\x\x c}{\smash[b]{\Z_{30}}}$ & $\Z_{120}\underset{c}{\oplus}\Z_2$ & $\underset{\x\x\x\x c}{\smash[b]{\Z_{120}}}$ & \multicolumn{1}{|c|}{$0$}\\
\hline
\end{tabular}
\end{center}
\medskip
\caption{Table of the first homotopy groups $\pi_r\U(n)$}
\end{table}

\bigskip

%
%

%

\end{document}